\documentclass[11pt]{amsart}
\usepackage{amsmath,amssymb,amscd,fullpage}

\theoremstyle{plain}

\newtheorem*{theorem*}{Theorem A}
\newtheorem*{theoremm*}{Theorem B}

\newtheorem{cor}[equation]{Corollary}

\newtheorem{rmks}[equation]{Remarks}

\numberwithin{equation}{subsection}

\newcommand{\C}{\mathbb{C}}

\newcommand{\A}{\mathbb{A}}

\newcommand{\N}{\mathbb{N}}

\newcommand{\Gal}{\operatorname{Gal}}

\title{Descent Construction For GSpin Groups: \\ Main Results and Applications}
\author{Joseph Hundley and Eitan Sayag}

\begin{document}
\maketitle

The purpose of this note is to announce an extension of the
descent method of Ginzburg, Rallis, and Soudry to the setting of
{\it essentially} self dual representations. This extension of the
descent construction provides a complement to recent work of
Asgari and Shahidi \cite{MR2219256}
on the generic transfer for general Spin groups as well as to the
work of Asgari and Raghuram \cite{AsgariRaghuram} on cuspidality
of the exterior square lift for representations of $GL_4$.
Complete proofs of the results announced in the present note will
appear in our forthcoming article(s).


\section{Preliminaries}
\subsection{GSpin groups and their quasisplit forms}
Let $F$ be a number field.  By the classification results in
Chapter 16 of \cite{MR1642713}
 (see especially 16.3.2,
16.3.3, 16.4.2), and the definition of the $L$ group, there is a
unique quasisplit $F$ group $G$ such that the connected component
of the identity in $^LG$ is $GSp_{2n}(\C).$  This is
$GSpin_{2n+1}.$

Similarly, there is a 1-1 correspondence between quasisplit $F$
groups $G$ such that $^LG^0 =GSO_{2n}(\C)$ and homomorphisms from
$\Gal( \bar F/ F)$ to the group with two elements, and hence, by
class field theory, with quadratic characters of $\chi:
\A_F^\times / F^\times \to \{\pm 1\}$ (the case $n=4$ is no
different, see Section \ref{finalrmrks}). The unique split group
$G$ such that $^LG^0 =GSO_{2n}(\C)$ corresponds to the trivial
character. We denote this group $GSpin_{2n}.$ The finite Galois
form of its $L$ group is $GSO_{2n}(\C).$ The form corresponding to
the nontrivial character $\chi$ we denote by $GSpin_{2n}^\chi.$
The finite Galois form of its $L$ group is $GSO_{2n}(\C) \rtimes
Gal( E/F)$ where $E$ is the quadratic extension of $F$
corresponding to $\chi.$

\subsection{Liftings}

According to the Langlands functoriality conjecture, one expects a
lifting of automorphic representations of $GSpin_{2n+1}( \A)$ to
automorphic representations of $GL_{2n}( \A)$ corresponding to the
inclusion
$$GSp_{2n}(\C) \hookrightarrow GL_{2n}(\C).$$
Similarly, one expects a lifting of automorphic representations of
$GSpin_{2n}( \A)$ to automorphic representations of $GL_{2n}( \A)$
corresponding to the inclusion
$$GSO_{2n}(\C) \hookrightarrow GL_{2n}(\C).$$
For (globally) generic representations, the existence of these
liftings is proved in \cite{MR2219256}.

Now consider $GSpin_{2n}^\chi$ for $\chi \ne 1.$
Regardless of $\chi,$ the $L$ group
$GSO_{2n}(\C) \rtimes Gal( E/F),$
with $E$ as above,
 is isomorphic
to $GO_{2n}(\C),$ and a specific isomorphism can be fixed by
mapping the nontrivial element of $Gal(E/F)$ to
$$\begin{pmatrix} I_{n-1} &&& \\ &&1& \\ &1&& \\ &&& I_{n-1}\end{pmatrix}.$$
One then expects a lifting of automorphic representations of
$GSpin_{2n}^\chi( \A)$ to automorphic representations of $GL_{2n}(
\A)$ corresponding to the inclusion
$$GO_{2n}(\C) \hookrightarrow GL_{2n}( \C).$$

We understand that this is the subject of \cite{A-S2}. We denote
all of these liftings by $AS.$

\section{Main results}
\subsection{The odd case}

\begin{theorem*}\label{main odd}Let $\omega$ be a Hecke character.
Suppose $n_1, \dots, n_m \in \N,$ and that, for each $1 \leq i
\leq m,$ $\tau_i$ is an irreducible cuspidal automorphic
representation of $GL_{2n_i}( \A)$ such that $L^S( s , \tau_i ,
\wedge^2 \times \omega^{-1} )$ has a pole at $s=1.$
Suppose furthermore that the representations $\tau_i$ are all distinct.
 Let $n=n_1+
\dots +n_m.$  Then there exists a globally generic irreducible
cuspidal automorphic representation $\sigma$ of $GSpin_{2n+1}(\A)$
such that $AS( \sigma) = Ind_{P_{\underline{n}}(\A)}^{GL_{2n}(\A)}
(\tau_1 \otimes \dots \otimes \tau_m )$ (normalized induction),
where $P_{\underline{n}}$ is the standard parabolic of $GL_{2n}$
corresponding to the ordered partition $2n= 2n_1+ \dots + 2n_m$ of
$2n.$ Furthermore, the central character of $\sigma$ is $\omega$.
\end{theorem*}

\begin{rmks}
\begin{enumerate}
\item
The notation ``$AS(\sigma)=\dots$'' requires some justification:
Theorem 1.1 of  \cite{MR2219256}
 assures the
{\it existence} of a weak lift $\Pi$ of $\sigma,$ but not its
uniqueness. However, for $\tau_1, \dots, \tau_m$ as in our
theorem, the induced representation
$Ind_{P_{\underline{k}}(\A)}^{GL_{2n}(\A)} (\tau_1 \otimes \dots
\otimes \tau_m ),$ is irreducible.  In conjunction with
Proposition 7.4, of \cite{MR2219256}
 this
implies that for the representation $\sigma$ which we construct,
 the transfer is uniquely determined.
\item
Note that if $\tau$ is a representation of $GL_\ell(\A)$ and $L^S(
s, \tau , \wedge^2 \times \omega^{-1})$ has a pole at $s=1$ then
it follows from \cite{MR1044830}
, Theorem 2, p. 224,
 that $\ell$ is even, and from  \cite{MR1044830}
 , Theorem 1, p. 213,
 that $\omega_{\tau}=\omega^{\frac{\ell}{2}}.$
\item
If $\tau$ is a representation of $GL_\ell(\A)$ such that $L^S(s,\tau, sym^2 \times
\omega^{-1})$ has a pole at $s=1,$ one may {\sl not} deduce that $\ell$ is
even.  However, one may deduce that $\tau \cong \tilde \tau \otimes \omega,$
whence $\omega^\ell = \omega_\tau^2$
(where $\omega_\tau$ is the central character of $\tau$).
  If $\ell$ is odd, it then follows that $\omega$
is the square of another global character $\eta,$ and that
$\tau'=\tau \otimes \eta^{-1}$ is self dual, with $L^S(s, \tau',
sym^2)$ having a pole at $s=1.$ Thus, the case when $\omega$ is a
square reduces to the self-dual case, and in the case when
$\omega$ is not a square we {\it can} deduce that $\ell$ is even
and that $\omega_{\tau}/{\omega^{\frac{\ell}{2}}}$ is quadratic.
\end{enumerate}
\end{rmks}

\subsection{The even case}
For the statement of the next main result, it will be convenient to
define $GSpin_{2n}^\chi:=GSpin_{2n}$ when $\chi$ is the trivial character.
\begin{theoremm*}\label{main even}
Let $\omega$ be a Hecke character which is not the square of
another Hecke character. Suppose $n_1, \dots, n_m \in \N,$ and
that, for each $i,$ $\tau_i$ is an irreducible cuspidal
automorphic representation of $GL_{2n_i}( \A)$ such that $L^S( s ,
\tau_i , sym^2 \times \omega^{-1} )$ has a pole at $s=1.$
Suppose furthermore that the representations $\tau_i$ are all distinct.
 Let
$n=n_1+ \dots +n_m,$ and, for each $i,$ let
$\chi_i=\omega_{\tau_i}/\omega^{n_i},$ which is quadratic.  Let
$\chi =\prod_{i=1}^m \chi_i.$    Then there exists a globally
generic irreducible cuspidal automorphic representation $\sigma$
of $GSpin^\chi_{2n}(\A)$ such that $AS( \sigma) =
Ind_{P_{\underline{n}}(\A)}^{GL_{2n}(\A)} (\tau_1 \otimes \dots
\otimes \tau_m )$ (normalized induction), where
$P_{\underline{n}}$ is the standard parabolic of $GL_{2n}$
corresponding to the ordered partition $2n= 2n_1+ \dots + 2n_m$ of
$2n.$ Furthermore, the central character of $\sigma$ is $\omega$.
\end{theoremm*}

\section{Applications}

\subsection{The image of the weak lift $AS$}

We now concentrate on the case of split general Spin groups. In
\cite{MR2219256}
, the authors show the existence
of functorial lifts from automorphic representations of
$GSpin_{2n}( \A)$ or $GSpin_{2n+1}( \A)$ to $GL_{2n}( \A).$ They
show that the images consist of automorphic representations which
satisfy the essential self-duality condition at almost all places.

Based on the self-dual case, (cf. Theorem A of \cite{MR1846354})
 one expects that the image of each
Asgari-Shahidi lifting consists of isobaric sums of distinct
essentially self dual cuspidal representations satisfying the
appropriate $L$-function condition. For example, any
representation in the image of the lift from $GSpin_{2n+1}$ should
be an isobaric sum of distinct $\omega$-symplectic cuspidals, for
some Hecke character $\omega$.


Our results support this expectation. We provide a ``lower bound''
for the image of the Asgari-Shahidi lifting, by showing that any
isobaric sum of distinct essentially self dual cuspidal
representations satisfying the appropriate $L$-function condition is
in the image of the appropriate lift.

\subsection{The image of the exterior square lift: $GL_4$ to $GL_6$}

The existence of an exterior square lift of a cuspidal automorphic
representation of $GL_4( \A)$ as an automorphic representation of
$GL_6( \A)$ was established by Kim in \cite{MR1937203}.
Recently, Asgari-Raghuram provided an explicit description of
those cuspidal automorphic representations of $GL_4( \A)$ whose
exterior square lift to $GL_6( \A)$ is not cuspidal. Among other
things their argument requires the following special case of
Theorem B.

\begin{cor}
Let $\Pi$ be a cuspidal representation of $GL_4( \A)$ and let
$\omega$ be any character of $GL_1(F)\backslash GL_1(\A)$. Assume that the partial
$L$-function $L^{S}(s,\Pi, sym^{2} \otimes \omega^{-1})$ has a pole
at $s = 1$ for a sufficiently large finite set $S$ of places of
$F$. Let $\chi=\omega_{\Pi}\omega^{-2}.$
Then there exists a globally generic cuspidal representation $\pi$
of $GSpin^{\chi}_4( \A)$ such that $\pi$ transfers to $\Pi$.
\end{cor}


Roughly speaking, Asgari and Raghuram prove that the exterior square lift of a
cuspidal representation $\Pi$ of $GL_4$ is cuspidal unless  $\Pi$
is isomorphic to a twist of either itself or its contragredient,
and that this occurs only if $\Pi$ is itself in the image of one
of four functorial lifts.  For the precise statement, see
\cite{AsgariRaghuram}, Theorem 1.1, p.2.  For the precise
relationship with our results, see p. 12.


\section{Scheme of Proof}

The proofs of Theorem {\bf A} and Theorem {\bf B} are obtained by
adapting (the special orthogonal group case of) the descent method
of Ginzburg, Rallis, and Soudry
\cite{MR1671452
,MR1722953
,MR1740991
,MR1846354
,MR1954940
}. The  adaptation is reasonably straightforward
owing to two observations:
\begin{enumerate}
\item There is a surjective homomorphism, defined over $F,$ from $GSpin_{m}$ to $SO_m,$
which restricts to an isomorphism between the unipotent
subvarieties.
\item The kernel of this projection is contained in the center of $GSpin_m.$  Thus, the
action $GSpin_m$ on itself by conjugation factors through the
projection.
\end{enumerate}

In what follows we detail the steps needed to prove Theorem {\bf
B}. The proof of Theorem {\bf A} is similar and technically
simpler.

The input to the construction is a collection $\underline \tau=
\left\{ \tau_1, \dots, \tau_m \right\}$ of cuspidal
representations $\tau_{i}$ of $GL_{2n_{i}}(\A)$ for $i=1, \dots, m,$
satisfying the assumptions of Theorem B.
Let $\chi_{\underline{\tau}}=\omega^{-n} \cdot \prod_{i=1}^m
\omega_{\tau_i}.$
Then $\chi_{\underline{\tau}}$ is a quadratic character.

 We can conveniently describe the
method in the following steps:

\begin{enumerate}
\item Construction of a family of descent representations 
of $GSpin^{\chi}_{4n+1-2\ell}(\A)$ for $\ell \geq n$.

\item Vanishing of the descent representations
for all $\ell > n$ and all $\chi \ne
 \chi_{\underline{\tau}}.$


\item Cuspidality and genericity (hence nonvanishing) of the descent representation of
$GSpin_{2n}^{\chi_{\underline \tau}}(\A)$.

\item Matching of spectral parameters at unramified places.
\end{enumerate}


The construction of the descent representations
relies on the notion of Fourier coefficient, as defined in
 \cite{MR1981592}, \cite{MR2214128}  (cf. also the ``Gelfand-Graev'' coefficients of
 \cite{MR2141707}).
    For purposes
 of presenting certain of the global arguments, it seems convenient to embed
 these Fourier coefficients into a slightly larger family of functionals,
 which we shall refer to as ``unipotent periods.''

  Suppose that $U$ is a unipotent subgroup of $GSpin_{4n+1}$ and
$\psi$ is a character of $U(F) \backslash U(\A).$  We define the
corresponding {\it unipotent period} to be the map from smooth,
left $U(F)$-invariant functions on $GSpin_{4n+1}(\A)$ to smooth,
left $(U(\A),\psi)$-equivariant functions,  given by
$$\varphi\mapsto \varphi^{(U,\psi)}$$
where
$$\varphi^{(U,\psi)}(g) := \int_{U(F)\backslash U(\A)} \varphi( ug ) \psi(u) \; du.$$

Each unipotent period has a local analogue at each finite place,
which is a twisted Jacquet functor.

Suppose now that $U$ is the unipotent radical of a standard
parabolic subgroup, and let $M$ denote the Levi.  The characters
of $U(F) \backslash U(\A)$ may be identified with the points of an
$F$-vector space, so that the stabilizer $Stab_M(\psi)$ makes
sense as an algebraic group defined over $F.$  We assume that
$\psi$ corresponds to a point in general position.
Then the map
$$FC^{\psi}: \left.\varphi \mapsto \varphi^{(U,\psi)}\right|_{Stab_M(\psi)(\A)}$$
is indeed a ``Fourier coefficient,'' as defined (and associated to
a nilpotent orbit) in \cite{MR1981592, MR2214128}.
  It maps smooth
functions of moderate growth on $GSpin_{4n+1}(F)\backslash
GSpin_{4n+1}(\A)$ to smooth functions of moderate growth on
$Stab_M(\psi)(F)\backslash Stab_M(\psi)(\A).$

Let $S$ be a set of
unipotent periods.  We will say that another unipotent period
$(U,\psi)$ {\it is spanned by} $S$ if
$$\left( \varphi^{(N,\vartheta)}\equiv 0 \; \forall (N,\vartheta) \in S \right)
\quad \Longrightarrow \quad \varphi^{(U,\psi)} \equiv 0.$$

We are now ready to describe each of the four steps listed  above in more
detail.

{\bf Step one: Construction of the descent
representations}

Using $\tau_1, \dots, \tau_m,$ a space of  Eisenstein series $E_{\underline
\tau,\omega}(g,\underline s)$ on $GSpin_{4n+1}(\A)$ is constructed-- corresponding
to a representation induced from  the standard parabolic subgroup $P=MU$ of
$GSpin_{4n+1}$ for which $M\cong GL_{2n_{1}} \times \dots \times
GL_{2n_{m}} \times GL_{1}.$
The partial $L$ functions
$$L^S( s, \tau_i, sym^2 \times \omega^{-1}) \qquad i=1 \text{ to } m$$
appear in the constant terms of elements of this space.
As a consequence, some of them have non-vanishing multi-residues at a
certain point $\underline s_0,$ {\it precisely because} of the $L$-function
hypothesis on $\underline \tau.$
In this fashion we obtain a residual representation-- which
lies in the discrete spectrum of $L^2( GSpin_{4n+1}(F)\backslash GSpin_{4n+1}(\A))$--
the nontriviality of which depends intrinsically on this $L$-function condition.
We denote  this representation by
$\mathcal{E}_{-1}(\underline{\tau},\omega)$.

Now, $GSpin_{4n+1}$ contains a family of parabolic subgroups
$Q_\ell = L_\ell N_\ell, \; \ell =1$ to $2n,$ with $L_{\ell}$
isomorphic to $GL_{1}^{\ell} \times G_{4n-2\ell+1},$ having the crucial
property that for each character $\psi$ of $N_\ell$ in general
position, the identity component of the group
$Stab_{L_\ell}(\psi)$ is isomorphic to one of the groups
$GSpin_{4n-2\ell}^\chi.$
Fixing specific isomorphisms, we may pull back  each Fourier coefficient
$$FC^\psi(\mathcal{E}_{-1}(\underline{\tau},\omega))$$
as described above,
to a space of functions defined on $GSpin_{4n-2\ell}^\chi(\A).$
There are many characters $\psi$ for a given value of $\ell$ and $\chi,$ but they comprise a
single orbit for the action of $L_\ell(F)$ by conjugation,
and the various spaces
$FC^\psi(\mathcal{E}_{-1}(\underline{\tau},\omega))$ all pull back to the same space of
functions on
$GSpin_{4n-2\ell}^\chi(\A),$  regardless of the choice of $\psi$ in
this orbit and regardless of the choice of isomorphism $GSpin_{4n-2\ell}^\chi\to Stab_{L_\ell}(\psi)^0.$
(For this, we require the extension of meromorphic continuation of Eisenstein series to non-$K$-finite
sections, provided in \cite{MR2402686}.)

In this manner we obtain a space of functions on $GSpin_{4n-2\ell}^\chi(\A)$ for each value of $\chi.$
The family of representations thus obtained comprises the descent representations.

{\bf Step two: Vanishing of other descents}

For $\ell > n,$ one shows that the above Fourier coefficients
vanish identically on our residue representation
$\mathcal{E}_{-1}( \underline\tau,\omega).$ The reason is local:
the corresponding twisted Jacquet module of the unramified constituent
of the corresponding local induced
representation vanishes.  The same is true if
$\ell = n,$ at any unramified place $v$ such that  the identity component of
$Stab_{L_n}(\psi)$ is not isomorphic to $GSpin_{2n}^\chi$ over $F_v.$

The remaining descent representation, corresponding to $\ell = n$ and $\chi=\chi_{\underline{\tau}},$
may now be referred to as ``the'' descent without ambiguity.



 {\bf Step Three: Cuspidality and genericity of the
descent}

Next we appeal to  global arguments which may be presented in
terms of ``identities of unipotent periods.''

Consider the unipotent
period  on $C^\infty(GSpin_{4n+1}(F) \backslash GSpin_{4n+1}(\A))$
which consists of taking the constant term with respect to the maximal parabolic
with Levi isomorphic to $GL_{2n} \times GL_1,$ and then taking a
Whittaker integral on the $GL_{2n}$ Levi.  It can be shown that
this unipotent period does not vanish on our residue representation
$\mathcal{E}_{-1}(\underline{\tau},\omega).$  One shows that this unipotent
period is, in fact, spanned by the periods corresponding to Whittaker integrals
on the descent representations  (as  $\ell \ge n$ and $\chi$ vary).

Having proved by local arguments that these periods
vanishes identically on the residue representation $\mathcal{E}_{-1}(\underline \tau,\omega),$
whenever $\ell >n$ or $\chi \ne \chi_{\underline\tau},$
we deduce that they do not vanish identically when $\ell =n$ and $\chi =\chi_{\underline\tau}.$
  This shows that the space $FC^\psi(\mathcal{E}_{-1}(\underline\tau,\omega))$ is not only nontrivial,
  but supports a nontrivial global Whittaker integral.

Next, consider the unipotent periods on $C^\infty(GSpin_{4n+1}(F)
\backslash GSpin_{4n+1}(\A))$ which consist of taking the
constant term with respect to the maximal parabolic with Levi
isomorphic to $GL_k \times GSpin_{4n-2k+1}$ for some $k,$ and
then, if $k$ is
even,  performing the integral one would use to define a descent
representation of $GSpin_{4n-2k+1},$ with some value of $\ell$
larger than $n- \frac k 2.$  Combining the vanishing results of
Step two with well-known facts from the theory of Eisenstein
series, we deduce that all of these periods vanish identically on
the residue representation $\mathcal{E}_{-1}(\underline
\tau,\omega).$  We then show that the unipotent period which
corresponds to taking the  constant term of one of the functions
in $FC^{\psi}(\mathcal{E}_{-1}(\underline{\tau},\omega))$ is in
their span.

It follows that all the functions in the descent are cuspidal.
At this point, we may deduce that the descent representation decomposes
discretely as a direct sum of irreducible cuspidal automorphic
representations, at least one of which is generic.  We select one such
component for the representation $\sigma$ of Theorem B.  What remains
is to show that $\sigma$ lifts weakly to
$Ind_{P_{\underline{n}}(\A)}^{GL_{2n}(\A)}
(\tau_{1} \otimes \dots \otimes \tau_{m} ).$

{\bf Step Four: Matching of spectral parameters at unramified
places}

For $\ell = n,$ at an unramified place, where the identity component of
$Stab_{L_n}(\psi)$ is isomorphic to $GSpin_{2n}^{\chi_{\underline\tau}},$
the twisted Jacquet module of the unramified constituent of
the local induced representation
 is isomorphic,  as a $Stab_{L_n}(\psi)(F_v)$-module
to a certain induced representation of $Stab_{L_n}(\psi)(F_v).$
When restricted to the identity component, this representation may
not be irreducible. Nevertheless, we are able to deduce
 that any nonzero
irreducible component of the Fourier coefficient must lift weakly to
$Ind_{P_{\underline{n}}(\A)}^{GL_{2n}(\A)}
(\tau_{1} \otimes \dots \otimes \tau_{m} ).$


\section{Final Remarks}\label{finalrmrks}

\begin{enumerate}
\item
When considering the identification of $GO_{2n}(\C)$ with
$GSO_{2n}(\C) \rtimes Gal( E/F),$ one could also map the
nontrivial element to
$$\begin{pmatrix} -I_{n-1} &&& \\ &&-1& \\ &-1&& \\ &&& -I_{n-1}\end{pmatrix}.$$
This produces a slightly different functorial lift corresponding
to the twist of the one we have chosen above by the quadratic
character $\chi.$  Theorem B is, of course, true for this
``alternate'' lifting, as well, since one may ``untwist.''
\item
These are essentially the only distinct extensions of the
inclusion
 $GSO_{2n}(\C)\hookrightarrow GL_{2n}(\C)$ to   $GSO_{2n}(\C)\rtimes Gal( E/F)$
 in the following sense.
 Suppose $V_1$ and $V_2$ are two $2n$ dimensional representations
 of $GSO_{2n}(\C)\rtimes Gal( E/F)$ such that the restriction of either to
 $GSO_{2n}(\C)$ is the standard representation.  Then
 one may show that $V_2$ is isomorphic to either $V_1$ or the
 twist of $V_1$ by the unique nontrivial character of $Gal(E/F).$
\item \label{d4remark}
A natural question arises in the case $n=4$: {\bf does the
$3$-fold symmetry of the Dynkin diagram of $GSpin_8$ lead to
additional quasi-split forms?} The answer is no, because the
$3$-fold symmetry of the $D_4$ root {\it system} does not extend
to a symmetry of the root {\it data} of $GSO_8$ and $GSpin_8.$
\end{enumerate}


\end{document}